\documentclass{amsart}

\usepackage{amsmath,amssymb,amsthm}
\usepackage{hyperref}
\usepackage{multicol}

\usepackage{graphicx}
\usepackage{xypic}
\entrymodifiers={+!!<0pt,\fontdimen22\textfont2>}
\usepackage[all]{xy}

\newtheoremstyle{myremark} 
    {7pt}                    
    {7pt}                    
    {}  	                 
    {}                           
    {\bf}       	         
    {.}                          
    {.5em}                       
    {}  

\theoremstyle{plain}
\newtheorem{lemma}{Lemma}
\newtheorem{theorem}[lemma]{Theorem}

\theoremstyle{myremark}

\newcommand{\longversion}[1]{}

\begin{document}
\title{The smallest nonevasive graph property}

\author[Micha{\l} Adamaszek]{Micha{\l} Adamaszek}
\address{Fachbereich Mathematik, Universit\"at Bremen
      \newline Bibliothekstr. 1, 28359 Bremen, Germany}
\email{aszek@mimuw.edu.pl}
\thanks{Author supported by a DFG grant.}



\begin{abstract}
A property of $n$-vertex graphs is called evasive if every algorithm testing this property by asking questions of the form ``is there an edge between vertices $u$ and $v$'' requires, in the worst case, to ask about all pairs of vertices. Most ``natural'' graph properties are either evasive or conjectured to be such, and of the few examples of nontrivial nonevasive properties scattered in the literature the smallest one has $n=6$.

 We exhibit a nontrivial, nonevasive property of $5$-vertex graphs and show that it is essentially the unique such with $n\leq 5$.
\end{abstract}
\maketitle

Evasiveness is a complexity-theoretic concept defined via the following combinatorial game. Two players, Alice and Bob, first fix a number $n$ and a property $\mathcal{P}$ of $n$-vertex graphs. Bob wants to find out if some unknown graph $G$, secretly chosen by Alice, has the property $\mathcal{P}$, by asking Alice one by one if a particular pair of vertices forms an edge. Alice wins if she can force Bob to ask about all the ${n\choose 2}$ pairs before he knows if $G\in\mathcal{P}$. Bob wins if he can decide the membership of $G$ in $\mathcal{P}$ after at most ${n\choose 2}-1$ questions. Of course there is no reason why Alice should fix any particular graph in advance --- she can adapt her answers so as to force Bob to ask the maximal number of questions. We say $\mathcal{P}$ is \emph{evasive} (or \emph{elusive}) if Alice has a winning strategy; it is \emph{nonevasive} if Bob does. For example, the simple property of ``being the complete graph'' is evasive. Alice's strategy is to say ``Yes'' to Bob's first ${n\choose 2}-1$ questions, at which point he is still not sure if $G$ is complete or not.

To be more precise, for a fixed natural number $n$ let $\mathcal{G}_n$ be the set of isomorphism classes of $n$-vertex simple, unlabeled graphs. A \emph{property} of $n$-vertex graphs is just an arbitrary subset $\mathcal{P}\subseteq \mathcal{G}_n$. We usually say ``a graph $G$ has property $\mathcal{P}$'' (e.g. $G$ is connected, $G$ is a tree, $G$ has a Hamiltonian cycle etc.) meaning ``$G$ is isomorphic to one of the graphs in $\mathcal{P}$''. For every $n$ there are two trivial nonevasive properties, $\mathcal{P}=\emptyset$ and $\mathcal{P}=\mathcal{G}_n$, for which Bob wins without asking any questions at all. More generally, $\mathcal{P}$ is evasive if and only if so is $\mathcal{G}_n\setminus\mathcal{P}$, with Bob playing the same strategy.

Evasiveness is a classical notion which arose as a way of measuring the decision-tree complexity of boolean functions. The lecture notes \cite{Lov} are an excellent introduction to this general topic. Here it suffices to say that most ``natural'' graph properties, for example connectedness, planarity, triangle-freeness, perfectness, existence of an isolated vertex and many more are all evasive. A major conjecture, attributed to Karp, claims that every nontrivial monotone property, that is a property closed under inserting new edges, is evasive. Its proof when $n$ is a prime power \cite{KKS} is one of the celebrated applications of topological methods in combinatorics.

Unsurprisingly, the known constructions of nonevasive properties are rare and to some extent artificial (see \cite{BEBL,MW76} for the original papers and \cite[Chapter 3]{Mark}, \cite[Chapter 13]{Koz} for surveys). The example usually presented in the literature involves classes of graphs called \emph{scorpions}, for which Bob can determine the answer after at most $6n-13$ questions, which is better than ${n\choose 2}$ for any $n\geq 11$. An optimized example of similar kind can be found in \cite[Fig.3.10]{Mark}. It is the property
\begin{center}
\raisebox{0.14cm}{$\mathcal{S}=\left\{\rule{0cm}{0.4cm}\right.\quad$}\includegraphics{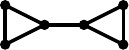}$\ $,$\qquad$\includegraphics{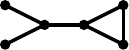}$\ $,$\qquad$\includegraphics{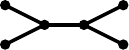}\raisebox{0.14cm}{$\quad\left\}\ \subseteq\ \mathcal{G}_6\rule{0cm}{0.4cm}\right.$}
\end{center}
of $6$-vertex graphs. Bob's strategy is to find the degree of one vertex, see where it could match the graphs in $\mathcal{S}$, and suitably expand his knowledge from there. It is an exercise to check that after $14$ out of the possible $15$ questions Bob can determine if Alice's graph has the desired shape except for not knowing the existence of one of the vertical edges. This, however, does not affect the membership in $\mathcal{S}$.

Clearly every nontrivial property of $3$-vertex graphs is evasive. Indeed, let $\mathcal{G}_3=\{G_0,G_1,G_2,G_3\}$ where $G_i$ is the unique $3$-vertex graph with $i$ edges. If $\mathcal{P}\subseteq\mathcal{G}_3$ is nontrivial then for some $i=0,1,2$ we have $G_i\in\mathcal{P}$ and $G_{i+1}\not\in\mathcal{P}$ or vice versa. Then Alice's strategy is to say ``Yes'' to Bob's first $i$ questions and ``No'' to the remaining ones. 

Is $n=6$ the smallest number of vertices for which there exists a nonevasive property? The answer may come as a surprise. It turns out that there is an essentially unique nonevasive property among graphs with at most $5$ vertices.

\begin{theorem}
\label{thm1}
The following property $\mathcal{E}\subseteq \mathcal{G}_5$ is nonevasive:
\begin{center}
\raisebox{0.37cm}{$\mathcal{E}=\left\{\rule{0cm}{0.6cm}\right.\quad$}\includegraphics{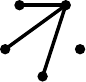},$\qquad$\includegraphics{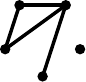},$\qquad$\includegraphics{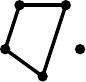},$\qquad$\includegraphics{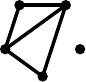},$\qquad$\includegraphics{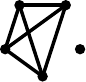},
\newline
\vskip0cm
$\qquad\qquad\qquad\quad$\includegraphics{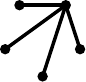},$\qquad$\includegraphics{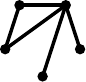},$\qquad$\includegraphics{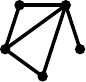},$\qquad$\includegraphics{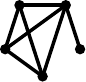},$\qquad$\includegraphics{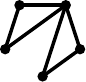},$\qquad$\includegraphics{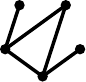}\raisebox{0.37cm}{$\quad\left\}\rule{0cm}{0.6cm}\right.$}.
\end{center}
Moreover, every nontrivial property of $4$-vertex graphs is evasive, while $\mathcal{E}$ and $\mathcal{G}_5\setminus \mathcal{E}$ are the only nontrivial, nonevasive properties of $5$-vertex graphs.
\end{theorem}
We chose $\mathcal{E}$ to stand for ``Exceptional'' or ``Eleven'', the cardinality of $\mathcal{E}$. 

\medskip
Let us sketch the proof of the theorem. A \emph{position} in the game is the complete graph on $n$ vertices whose edges are labeled with either ``present'', ``absent'' or ``unknown''. The first two indicate the status of an edge already discovered by Bob. The edges labeled ``unknown'' are those Bob hasn't asked about yet. A position with just one unknown edge is winning for Bob if the two graphs obtained by declaring this edge present or absent are either both in $\mathcal{P}$ or both not in $\mathcal{P}$. To find the winning player and winning moves for other positions we use the standard algorithm processing the game tree bottom-up. Now all is just a matter of efficiency. For $n=5$ there are $34$ graphs and $758$ isomorphism classes of positions. Evaluating the initial position, with all edges unknown, against all $2^{|\mathcal{G}_5|}=2^{34}$ graph properties, is a matter of at most one day on any reasonably modern personal computer. For $n=4$ the same thing is immediate. That verifies the theorem.

\medskip
We will now present Bob's winning strategy for $\mathcal{E}$ in a more accessible form. The property $\mathcal{E}$ has a special feature which reduces the number of positions we have to consider. For a graph $G$ let $\overline{G}$ denote its complement, that is the graph whose edges are the non-edges of $G$. If $P$ is a position in the game, we define the complement $\overline{P}$ by renaming all edges labeled ``present'' to ``absent'' and vice-versa. The edges unknown in $P$ remain unknown in $\overline{P}$.  Note that $\mathcal{E}$ contains five pairs of complements and one self-complementary graph. In other words $G\in\mathcal{E}$ if and only if $\overline{G}\in\mathcal{E}$. An easy inductive argument implies that in the $\mathcal{E}$-game a position $P$ is winning for Bob if and only if its complement $\overline{P}$ is. It is also easy to read off the strategy for $P$ from the strategy for $\overline{P}$. It follows that in our analysis we can identify a position with its complement. We can, for example, choose to work only with positions that have at least as many ``present'' as ``absent'' edges.

\begin{center}
\begin{figure}[h]
\includegraphics[scale=0.7]{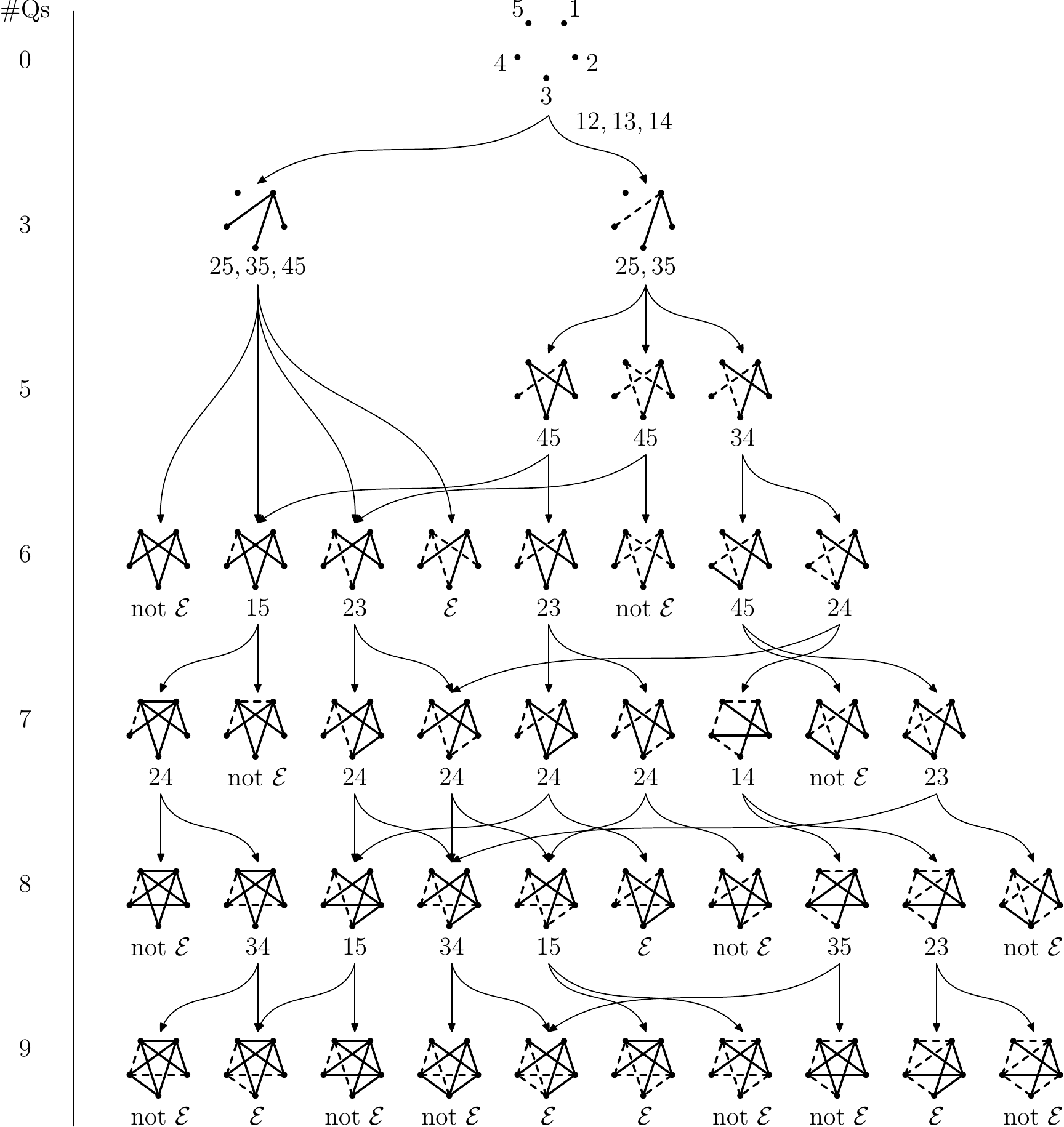}
\caption{Bob's winning strategy.\label{figure:bob}}
\end{figure}
\end{center}

The strategy is depicted in Figure~\ref{figure:bob} and here we offer a short description of how it starts. For convenience we label the vertices $1,\ldots,5$ as in the top of the figure. In the first three questions Bob asks about the edges $12$, $13$ and $14$. Up to isomorphism and complementation that leads to one of the positions in the second row of the figure. Depending on the outcome Bob now asks about three edges around vertex $5$ or two edges around $5$ and one other edge. After six questions we reach one of $8$ possible positions (again up to complementation), of which three already end the game. For the remaining three questions the rules are more complicated and it is best just to follow the arrows, each of which corresponds to one of Alice's answers, remembering that one may have to apply isomorphism and complementation of positions along the way. In each position the solid lines denote the present edges and the dashed lines are the absent edges (the edges not shown are unknown). The label under a position is either the next question to ask or an indication that membership in $\mathcal{E}$ is already decided. In each case Bob wins after the $9$th question at the latest.

\subsection*{Acknowledgement} The author thanks the Center for Mathematical Culture and the organizers of the $50$th Szkola Matematyki Pogladowej (\texttt{www.msn.ap.siedlce.pl}) for their invitation to give a talk whose preparation prompted the search for small examples of this kind.


\end{document}